\documentclass{article}
\usepackage{amsfonts}
\usepackage{amsmath}

\newtheorem{theorem}{Theorem}

\newtheorem{corollary}[theorem]{Corollary}

\newtheorem{definition}[theorem]{Definition}
\newtheorem{example}[theorem]{Example}

\newtheorem{lemma}[theorem]{Lemma}

\newtheorem{proposition}[theorem]{Proposition}
\newtheorem{remark}[theorem]{Remark}

\newenvironment{proof}[1][Proof]{\noindent\textbf{#1.} }{\ \rule{0.5em}{0.5em}}
\input{tcilatex}

\begin{document}

\title{An equivalent condition for a uniform space to be coverable}
\author{Conrad Plaut \\
Mathematics Department\\
University of Tennessee\\
Ayres Hall 121\\
Knoxville, TN 37996-1300\\
cplaut@math.utk.edu}
\maketitle

\begin{abstract}
We prove that an equivalent condition for a uniform space to be coverable is
that the images of the natural projections in the fundamental inverse system
are uniformly open in a certain sense. As corollaries we (1) obtain a
concrete way to find covering entourage, (2) correct an error in \cite{BPUU}%
, and (3) show that coverable is equivalent to chain connected and uniformly
joinable in the sense of \cite{RC}. Keywords: universal cover, uniform
space, coverable, fundamental group MSC: 55Q52; 54E15,55M10
\end{abstract}

In \cite{P} and \cite{BPUU} we approached the problem of defining universal
covers of locally bad spaces using the following ideas: (1) The appropriate
category in which to work in is not topological spaces, but rather uniform
spaces and uniformly continuous mappings. (2) The replacement for quotient
mappings in this category is bi-uniformly continuous mappings (see below).
(3) The appropriate replacement for curves is equivalence classes of chains.
We showed in \cite{BPUU} that such a program can be carried out for a large
class of uniform spaces called \textit{coverable spaces. }In particular we
constructed, for any uniform space $X$, a uniform space $\widetilde{X}$, a
natural uniformly continuous mapping $\phi :\widetilde{X}\rightarrow X$ and
a group $\delta _{1}(X)$ that acts on $\widetilde{X}$. For coverable spaces
the mapping $\phi $ has many of the properties of a universal covering map,
such as lifting and universal properties, and we refer to the space $%
\widetilde{X}$ as the uniform universal cover of $X$. The group $\delta
_{1}(X)$ (which we called the \textquotedblleft deck
group\textquotedblright\ in \cite{BPUU} but which was renamed the
\textquotedblleft uniform fundamental group\textquotedblright\ in \cite{PL})
is a functorial invariant of uniform structures having properties like the
fundamental group in this category. See \cite{BPUU} for many specific
theorems and examples, and \cite{BPCR} for some additional applications.

The space $\widetilde{X}$ is the inverse limit of the fundamental inverse
system $(X_{E},\phi _{EF})$ of $X$, which is indexed on the set of all
entourages $E$ of $X$. Roughly speaking, $X_{E}$ \textquotedblleft
unrolls\textquotedblright\ nontrivial classes of loops that are in some
sense larger than $E$ (we will give more background below). We will denote
by $\phi _{E}:\widetilde{X}\rightarrow X_{E}$ and $\phi :\widetilde{X}%
\rightarrow X$ the natural projections (the latter is actually just the
endpoint mapping). $X$ is coverable by definition if $\phi $ is surjective $%
\phi _{E}$ is surjective for all $E$ in some basis (called a covering basis)
for the uniform structure of $X$. Elements of the covering basis are called
covering entourages. Existence of a covering basis can be proved in many
cases--for example it is not too difficult to show that connected and
uniformly locally connected pseudometric spaces--which includes all geodesic
spaces--are coverable (Theorem 98, \cite{BPUU}). On the other hand, it is
easy to find in coverable spaces entourages that are not covering entourages
(cf. Example \ref{exex}), and without uniform nice local conditions it can
be difficult to verify coverability. In this paper we show that coverability
is equivalent to the following: $X$ is chain connected and for any entourage 
$E$, $\phi _{E}$ has image that is \textit{uniformly open }in $X_{E}$ in a
sense defined below (Theorem \ref{main}). While surjectivity of maps in an
inverse system is generally a strong condition useful for proving theorems,
from the standpoint of verification there is a clear advantage in not having
to hunt for the covering entourages. Theorem \ref{main} provides a condition
that may be verified for \textit{arbitrary }entourages. Moreover, as a
corollary we obtain a constructive method for extracting a covering
entourage from an arbitrary entourage without having to consider the mapping 
$\phi _{E}$ at all (see Corollary \ref{extract} and Example \ref{exex}).

In \cite{RC} the authors explore our construction of the uniform universal
cover in the setting of what they call \textit{uniformly joinable} uniform
spaces. At first their construction does not look like our construction, and
one must look in Section 8 of \cite{RC} to find a statement that they are
\textquotedblleft identical\textquotedblright . As we explain below, their
definitions of $GP(X,\ast )$ and $\overset{\vee }{\pi }(X,\ast )$ are simply
translations of our definitions of $\widetilde{X}$ and $\delta _{1}(X)$ into
the language of Rips complexes. Moreover, another application of our main
theorem is that the class of chain connected, uniformly joinable spaces
considered in \cite{RC} is precisely the same as the class of coverable
spaces (Corollary \ref{eqrc}). The authors of \cite{RC} state that
\textquotedblleft a topologist would be skeptical\textquotedblright\ of this
particular result--an assertion backed by some musings on Siebenmann's
thesis. Nonetheless, their definition is closely related to concepts in
continua theory and they obtain the interesting result that a metric
compactum $X$ is uniformly joinable if and only if the function $\phi :%
\widetilde{X}\rightarrow X$ is surjective. In light of Corollary \ref{eqrc},
for metrizable spaces this is a strong generalization of the fact, proved in 
\cite{BPLC}, that a compact topological group is coverable if and only if $%
\phi $ is surjective. Also in \cite{RC} the authors introduce a notion of
generalized cover in the uniform category that does not require a group
action. In light of this, what we called \textquotedblleft covers of uniform
spaces\textquotedblright\ in \cite{BPUU} really should be called something
like \textquotedblleft regular uniform covers\textquotedblright\ as is
suggested in \cite{RC}.

We do not use any theorems from \cite{RC} in our proofs and in fact, in
light of Corollary \ref{eqrc}, some of the theorems in \cite{RC} were
already proved in \cite{BPUU}. Thanks to Jurek Dydak for critiques and
stimulating comments. In particular, he pointed out an error in \cite{BPUU}
that is corrected in the present paper. Valera Berestovskii provided some
valuable comments.

We will use the notation of \cite{BPUU}. In particular, we generally use $f$
in place of $f\times f$; for example, if $E$ is an entourage in a uniform
space we will write $f(E)$ rather than $\left( f\times f\right) (E)$. For an
entourage $E$ we let $B(x,E):=\{y:(x,y)\in E\}$. Given a uniform space $X$,
for any entourage $E$, $X_{E}$ is defined to be the space of $E$-homotopy
classes of $E$-chains $\alpha :=\{x_{0}=\ast ,...,x_{n}\}$, where $\ast $ is
a basepoint. By definition, $\alpha $ is an $E$-chain if $(x_{i},x_{i+1})\in
E$ for all $i$. An $E$-homotopy of $\alpha $ is a finite sequence of moves
starting with $\alpha $, where each move consists of adding or taking away a
point (but not endpoints!) so long as doing so results in an $E$-chain. For
chain connected spaces (meaning every pair of points is joined by an $E$%
-chain for all $E$) nothing of consequence depends on the choice of
basepoint so we generally eliminate it from the notation. The space $X_{E}$,
the elements of which are denoted $[\alpha ]_{E}$, has a natural uniform
structure having a basis consisting of sets $F^{\ast }$, where $F\subset E$
and $([\alpha ]_{E},[\beta ]_{E})\in F^{\ast }$ if and only if $[\alpha
]_{E}=[x_{0},...,x_{n-1},x]_{E}$, and $[\beta
]_{E}=[x_{0},...,x_{n-1},y]_{E} $, with $(x,y)\in F$. $\widetilde{X}$ is
given the inverse limit uniformity. When $F\subset E$, the mapping $\phi
_{EF}:X_{F}\rightarrow X_{E}$ simply considers an $F$-chain as an $E$-chain,
i.e., $\phi _{EF}([\alpha ]_{F})=[\alpha ]_{E}$, and $\phi
_{XE}:X_{E}\rightarrow X$ is the endpoint mapping. With respect to the
natural uniform structure these mappings are bi-uniformly continuous in the
sense that the inverse image of any entourage is an entourage, and the image
of any entourage is an entourage in the subspace uniformity of the image of
the mapping. Given a uniformly continuous mapping $f:X\rightarrow Y$ and
entourages $E,F$ in $X,Y$, respectively, such that $f(E)\subset F$, there is
a unique basepoint-preserving induced uniformly continuous function $%
f_{EF}:X_{E}\rightarrow Y_{F}$ such that $\phi _{YF}\circ f_{EF}=f\circ \phi
_{XE}$, which simply takes $[\alpha ]_{E}$ to $[f(\alpha )]_{F}$.

If $X$ is chain connected, the function $\phi _{XE}:X_{E}\rightarrow X$ is a
quotient mapping via the isomorphic action of the group $\delta _{E}(X)$
consisting of $E$-homotopy classes of $E$-loops. Precisely what this means
is not needed for this paper (see \cite{P} for the definitions); we can get
by with two facts: first, if $\phi _{XE}(a)=\phi _{XE}(b)$ then for some $%
g\in \delta _{E}(X)$, $g(a)=b$ and second, the entourages $F^{\ast }$ are 
\textit{invariant} in the sense that for every $g\in \delta _{E}(X)$, $%
g(F^{\ast })=F^{\ast }$.

\begin{definition}
We say that a subset $A$ of a uniform space $X$ is uniformly open if there
is an entourage $E$ in $X$ such that for every $a\in A$, $B(a,E)\subset A$.
\end{definition}

There are a few obvious facts: if $A$ is a uniformly open set then $A$ is
open, the complement of $A$ is uniformly open, and hence $A$ itself is also
closed. But for example in the rational numbers $\mathbb{Q}$ with the usual
metric there are plenty of open and closed subsets that are not uniformly
open. The inverse image of any uniformly open set via a uniformly continuous
function is uniformly open, but in general nothing can be said of images.
For example consider the bi-uniformly continuous surjection $f:[0,2]\mathbb{%
\times Z}_{2}\rightarrow \lbrack 0,2]$ defined by $\mathbb{(}q,0)\longmapsto
q$ and $(q,1)\longmapsto \frac{q}{2}$. Here $[0,2]$ has its usual metric, $%
\mathbb{Z}_{2}$ has the discrete metric, and $[0,2]\mathbb{\times Z}%
_{2}\rightarrow \lbrack 0,2]$ has the product metric. It is easy to check
that $[0,2]\times \{1\}$ is uniformly open in $[0,2]\mathbb{\times Z}_{2}$
but of course $f([0,2]\times \{1\})=[0,1]$ is not even open in $[0,2]$. (But
see Remark \ref{quot} below.)

\begin{lemma}
\label{uolem}A uniform space $X$ is chain connected if and only if the only
non-empty uniformly open subset of $X$ is $X$.
\end{lemma}

\begin{proof}
Suppose $X$ is chain connected and let $U$ be a nonempty uniformly open
subset of $X$. If $E$ is an entourage as in the definition of uniformly
open, then any $E$-chain starting at $x$ cannot leave $U$ and so $U=X$. For
the converse, let $x\in X$, $E$ be an entourage, and $U$ be the set of all
points that are joined to $x$ by an $E$-chain. If $z\in U$ then clearly $%
B(z,E)\in U$; hence $U$ is uniformly open and non-empty, hence equal to $X$.
Since $E$ and $x$ were arbitrary, $X$ is chain connected.
\end{proof}

Obviously the intersection of any two uniformly open subsets is uniformly
open. As a corollary of this and the above lemma we obtain:

\begin{corollary}
\label{cool}If any two uniformly open chain connected subsets of a uniform
space $X$ have non-empty intersection then they must be equal.
\end{corollary}

\begin{definition}
If $X$ is a uniform space, $E\subset F$ are entourages in $X$, and $A$ is a
uniformly open subset of $X_{E}$, define 
\begin{equation*}
F_{A}:=\phi _{XE}(F^{\ast }\cap \left( A\times A\right) )\text{.}
\end{equation*}
\end{definition}

\begin{lemma}
\label{greatlem}Let $X$ be a chain connected uniform space, $E$ be an
entourage in $X$, and $A$ be a uniformly open subset of $X_{E}$. Then $\phi
_{XE}(A)=X$ and for any entourage $F\subset E$, $F_{A}$ is an entourage in $%
X $.
\end{lemma}

\begin{proof}
Consider entourages $W\subset F\subset E$ such that if $x\in A$ and $%
(x,y)\in W^{\ast }$ then $y\in A$. We will first prove that $\phi _{XE}(A)$
is uniformly open and hence equal to $X$. Suppose that $(a,b)\in W$ and $%
a\in \phi _{XE}(A)$. So there exist $z\in A$ such that $\phi _{XE}(z)=a$ and
since $W=\phi _{XE}(W^{\ast })$ (cf. \cite{BPUU}, Proposition 16) there
exists $(x,y)\in W^{\ast }$ such that $\phi _{XE}(x,y)=(a,b)$. Next there
exists some $g\in \delta _{E}(X)$ such that $g(x)=z$. By the invariance of $%
W^{\ast }$, if $w:=g(y)$ then $(z,w)\in W^{\ast }$ and $\phi
_{XE}(z,w)=(a,b) $. By choice of $W^{\ast }$, $w\in A$, which places $b\in
\phi _{XE}(A)$, finishing the proof that $\phi _{XE}(A)$ is uniformly open
and equal to $X$. Now the initial assumption that $a\in \phi _{XE}(A)$ is
superfluous and the same argument shows $W\subset W_{A}$ and since $%
W_{A}\subset F_{A}$, $F_{A}$ is an entourage.
\end{proof}

\begin{remark}
\label{quot}The same proof as in the previous lemma shows the following: If $%
f:X\rightarrow Y$ is a quotient of uniform spaces via an isomorphic action
and $A\subset X$ is uniformly open then $f(A)$ is uniformly open in $Y$ (see 
\cite{P} for a discussion of isomorphic actions).
\end{remark}

\begin{corollary}
\label{greatcor}If $X$ is a chain connected uniform space and there is some
entourage $E$ such that $\phi _{E}(\widetilde{X})$ is uniformly open in $%
X_{E}$ then $\phi :\widetilde{X}\rightarrow X$ is surjective.
\end{corollary}

\begin{lemma}
\label{form}Let $X$ be a chain connected uniform space and $E$ be an
entourage in $X$ such that $A:=\phi _{E}(\widetilde{X})$ is uniformly open.
If $\alpha :=\{\ast =x_{0},...,x_{n}\}$ is an $E_{A}$-chain then $[\alpha
]_{E}\in A$.
\end{lemma}

\begin{proof}
We will show by induction that $[x_{0},...,x_{k}]_{E}\in A$ for all $k\leq n$%
. Certainly the statement is true for $k=0$. Suppose that $%
[x_{0},...,x_{k}]_{E}\in A$. Now $(x_{k},x_{k+1})\in E_{A}$ and by
definition there exist $E$-chains $\gamma :=\{\ast =y_{0},...,y_{m},x_{k}\}$
and $\omega :=\{y_{0},...,y_{m},x_{k+1}\}$ such that $[\gamma ]_{E},[\omega
]_{E}\in A$ and $(x_{k},x_{k+1})\in E$. So we have $\left(
[a_{D}]_{D}\right) ,([b_{D}]_{D}),([c_{D}]_{D})\in \widetilde{X}$ such that $%
[a_{E}]_{E}=[x_{0},...,x_{k}]_{E}$, $[b_{E}]_{E}=[\gamma ]_{E}$, and $%
[c_{E}]_{E}=[\omega ]_{E}$. Consider $\kappa :=([a_{D}\ast b_{D}^{-1}\ast
c_{D}]_{D})\in \widetilde{X}$, where \textquotedblleft $\ast $%
\textquotedblright\ denotes concatenation of chains. Now 
\begin{equation*}
\phi _{E}(\kappa )=[a_{E}\ast b_{E}^{-1}\ast
c_{E}]_{E}=[\{x_{0},...,x_{k}\}\ast \gamma ^{-1}\ast \omega ]_{E}
\end{equation*}%
\begin{equation*}
=[x_{0},...,x_{k},y_{m},...,y_{1},y_{0},y_{1},...,y_{m},x_{k+1}]_{E}=[x_{0},...,x_{k+1}]_{E}
\end{equation*}%
where the last $E$-homotopy successively removes $%
y_{0},y_{1},y_{1},y_{2},...,y_{m}$.
\end{proof}

\begin{proposition}
\label{diag}Let $X$ be a chain connected uniform space and $E$ be an
entourage in $X$ such that $A:=\phi _{E}(\widetilde{X})$ is uniformly open.
Letting $D:=E^{\ast }\cap (A\times A)$ and $G:=\phi _{E}^{-1}(D)=\phi
_{E}^{-1}(E^{\ast })$, there is a uniformly continuous function $\psi :%
\widetilde{X}\rightarrow A_{D}$ such that the following diagram commutes%
\begin{equation}
\begin{array}{lll}
\widetilde{X}_{G} & \overset{\phi _{\widetilde{X}G}}{\longrightarrow } & 
\widetilde{X} \\ 
\downarrow ^{\theta } & \swarrow _{\psi } & \downarrow ^{\phi _{E}} \\ 
A_{D} & \overset{\phi _{AD}}{\longrightarrow } & A%
\end{array}
\label{diagram}
\end{equation}%
where $\theta =(\phi _{E})_{GD}$ is the mapping induced by $\phi _{E}$.
\end{proposition}

\begin{proof}
First of all we recall what it means for $\left( [\alpha ]_{E},[\beta
]_{E}\right) \in D$: $[\alpha ]_{E},[\beta ]_{E}\in A$ and 
\begin{equation*}
\lbrack \alpha ]_{E}=[y_{0},...,y_{m},x]_{E}\text{ and }[\beta
]_{E}=[y_{0},...,y_{m},y]_{E}
\end{equation*}%
for some choice of $y_{0},...,y_{m}$, with $\left( x,y\right) \in E$. We
will define $\psi :=f\circ \phi _{E_{A}}$, where 
\begin{equation}
f([x_{0},...,x_{n}]_{E_{A}})=[[x_{0}]_{E},[x_{0},x_{1}]_{E},...,[x_{0},...,x_{n}]_{E}]_{D}
\label{map}
\end{equation}%
We need to check various things about $f$. First of all, note that by Lemma %
\ref{form}, $[x_{0},...,x_{k}]_{E}\in A$ for all $1\leq k\leq n$ and
therefore 
\begin{equation*}
([x_{0},...,x_{k}]_{E},[x_{0},...,x_{k},x_{k+1}]_{E})
\end{equation*}%
\begin{equation*}
=([x_{0},...,x_{k},x_{k}]_{E},[x_{0},...,x_{k},x_{k+1}]_{E})\in D
\end{equation*}%
so the definition at least goes into the correct set. To see that $f$ is
well-defined, consider an $E_{A}$-chain $\alpha ^{\prime
}:=\{x_{0},...,x_{k-1},x,x_{k},...,x_{n}\}$, which would lead to 
\begin{equation*}
\lbrack \lbrack
x_{0}]_{E},...,[x_{0},...,x_{k-1}]_{E},[x_{0},...,x_{k-1},x]_{E},
\end{equation*}%
\begin{equation*}
\lbrack
x_{0},...,x_{k-1},x,x_{k}]_{E},...,[x_{0},...,x_{k-1},x,x_{k},...,x_{n}]_{E}]_{D}
\end{equation*}%
in the above definition. But notice that for any $k\leq m\leq n$, we already
know $\{x_{0},...,x_{k-1},x,x_{k},...,x_{m}\}$ and $%
\{x_{0},...,x_{k-1},x_{k},...,x_{m}\}$ are $E_{A}$-homotopic, hence $E$%
-homotopic, so we can simply remove the \textquotedblleft $x$%
\textquotedblright\ from all such terms. This leaves the one extra term $%
[x_{0},...,x_{k-1},x]_{E}$. But since $%
([x_{0},...,x_{k-1}]_{E},[x_{0},...,x_{k-1},x_{k}]_{E})\in D$, up to $D$%
-homotopy we may simply remove this term, getting us back to (\ref{map}).

We will now check that $f$ is uniformly continuous. To do so we will have to
be a little more careful with notation. Given an entourage $W\subset E_{A}$
in $X$ we have an entourage called $W^{\ast }$ in $X_{E}$ and one called $%
W^{\ast }$ in $X_{E_{A}}$. We will refer to the latter as $W^{\#}$. We also
have the entourage $(W^{\ast }\cap (A\times A))^{\ast }$ in $A_{D}$, which
we will simply denote by $W^{\ast \ast }$. The proof of uniform continuity
will be finished if we can show that $f(W^{\#})\subset W^{\ast \ast }$. Let $%
([\alpha ]_{E_{A}},[\beta ]_{E_{A}})\in W^{\#}$. By definition we may take $%
\alpha =\{x_{0},...,x_{n},x\}$ and $\beta =\{x_{0},...,x_{n},y\}$ with $%
(x,y)\in W$. We have 
\begin{equation*}
f([\alpha
]_{E_{A}})=[[x_{0}]_{E},...,[x_{0},...,x_{n}]_{E},[x_{0},...,x_{n},x]_{E}]_{D}
\end{equation*}%
\begin{equation*}
f([\beta
]_{E_{A}})=[[x_{0}]_{E},...,[x_{0},...,x_{n}]_{E},[x_{0},...,x_{n},y]_{E}]_{D}
\end{equation*}%
Since $(x,y)\in W$ and $[x_{0},...,x_{n},x]_{E},[x_{0},...,x_{n},y]_{E}\in A$
(Lemma \ref{form} again), $([x_{0},...,x_{n},x]_{E},[x_{0},...,x_{n},y]_{E})%
\in W^{\ast }\cap (A\times A)$ and $\left( f([\alpha ]_{E_{A}}),f([\beta
]_{E_{A}})\right) \in W^{\ast \ast }$.

We will now check the commutativity of the diagram. Suppose that $\eta
:=\{y_{0}=\ast ,y_{1},...,y_{n}\}$ is a $G$-chain in $\widetilde{X}$. This
means that for all $i$, $\left( \phi _{E}(y_{i}),\phi _{E}(y_{i+1})\right)
\in D$. In particular, $\left( \phi _{E}(y_{0}),\phi _{E}(y_{1})\right)
=\left( [\ast ]_{E},\phi _{E}(y_{1})\right) \in D$. This means that we may
write $\phi _{E}(y_{1})=[\ast =w_{0},...,w_{m},x_{1}]_{E}$, where $x_{1}$ is
the endpoint of $y_{1}$, $\{w_{0},...,w_{m},\ast \}$ is $E$-homotopic to the
identity and $(\ast ,x_{1})\in E$. But then we may use the null $E$-homotopy
of $\{w_{0},...,w_{m},\ast \}$ to see that 
\begin{equation*}
\phi _{E}(y_{1})=[\ast ,...,w_{m},\ast ,x_{1}]_{E}=[\ast ,x_{1}]_{E}\text{.}
\end{equation*}%
By definition of $D$ we also have $[\ast ,x_{1}]_{E}=\phi _{E}(y_{1})\in A$,
which implies that $\{\ast ,x_{1}\}$ is an $E_{A}$-chain. Proceeding
inductively with essentially the same argument, we see that $\phi
_{E}(y_{i})=[\ast ,x_{1},..,x_{i}]_{E}$, where $x_{i}$ is the endpoint of $%
\phi _{E}(y_{i})$ and $\{x_{0},...,x_{n}\}$ is an $E_{A}$-chain. By
definition of $\theta $, 
\begin{equation*}
\theta ([\eta ]_{G})=[\phi _{E}(y_{0}),...,\phi _{E}(y_{n})]_{D}=[[\ast
]_{E},[\ast ,x_{1}]_{E},...,[\ast ,x_{1},...,x_{n}]_{E}]_{D}
\end{equation*}%
\begin{equation*}
=\psi (y_{n})=\psi \circ \phi _{\widetilde{X}G}([\eta ]_{G})
\end{equation*}%
This proves the commutativity of the upper triangle. The commutativity of
the lower triangle is obvious from the definition of $\psi $.
\end{proof}

Universal uniform spaces and universal bases were defined in \cite{BPUU};
the definitions will be explained in the proof below.

\begin{proposition}
\label{up}If $X$ is a chain connected uniform space such that for every
entourage $E$, $\phi _{E}:\widetilde{X}\rightarrow X_{E}$ has uniformly open
image in $X_{E}$ then $\widetilde{X}$ is universal with an invariant (with
respect to the action of $\delta _{1}(X)$) universal basis.
\end{proposition}

\begin{proof}
Consider the diagram (\ref{diagram}). We will start by showing that $%
\widetilde{X}$ is chain connected. Since $\phi _{E}$ is surjective onto $A$,
so is $\phi _{AD}$. This means that every pair of points in $A$ is joined by
a $D$-chain. Equivalently, $A\times A=\bigcup\limits_{n=1}^{\infty }D^{n}$
(where $D^{n}$ is the set of all points in $A$ joined to the basepoint by a $%
D$-chain of length $n$). Since $\phi _{E}$ is surjective, it is easy to
check that $\phi _{E}^{-1}(D^{n})=\left( \phi _{E}^{-1}(D)\right)
^{n}=\left( \phi _{E}^{-1}(E^{\ast })\right) ^{n}$ and so $\widetilde{X}%
=\bigcup\limits_{n=1}^{\infty }\left( \phi _{E}^{-1}(E^{\ast })\right) ^{n}$
(cf. the proof of Lemma 11 in \cite{BPUU}). This means that every pair of
points in $\widetilde{X}$ is joined by a $\phi _{E}^{-1}(E^{\ast })$-chain.
Since the set of all $\phi _{E}^{-1}(E^{\ast })$ forms a basis for the
uniformity of $\widetilde{X}$, $\widetilde{X}$ is chain connected. This now
implies that the mapping $\phi _{\widetilde{X}G}$ is surjective and the
hypotheses of Proposition 33 in \cite{BPUU} are satisfied for this diagram,
implying that $\phi _{\widetilde{X}G}$ is a uniform homeomorphism. By
definition $G$ is a universal entourage, and since it is of the form $\phi
_{E}^{-1}(E^{\ast })$, it is invariant (cf. Proposition 41 of \cite{P}). We
have shown that $\widetilde{X}$ has a basis of universal entourages; by
definition this makes $\widetilde{X}$ universal.
\end{proof}

Alas, the proof of Corollary 61 in \cite{BPUU} is not correct--or rather,
the proof is correct for a weaker statement. The penultimate sentence in the
proof requires an additional assumpion. For example, the proof is correct
for the following statement:

\begin{lemma}
\label{61}If $f:X\rightarrow Y$ is a quotient via an action on a uniform
space $X$ and $X$ has a universal basis that is invariant with respect to
the action, then $Y$ is coverable.
\end{lemma}

Corollary 61 was only used to establish the equivalence of the definition of
\textquotedblleft coverable topological group\textquotedblright\ as defined
in \cite{BPTG} with the definition in \cite{BPUU} when applied to
topological groups considered as uniform spaces. In particular, none of the
results of \cite{BPTG} cited in the current paper relies on this corollary;
Lemma \ref{61} will suffice to prove our main Theorem \ref{main}, of which
Corollary 61 is a corollary.

\begin{theorem}
\label{main}For a chain connected uniform space $X$, the following are
equivalent:

\begin{enumerate}
\item $X$ is coverable.

\item $\phi :\widetilde{X}\rightarrow X$ is a bi-uniformly continuous
surjection.

\item For each entourage $E$ in $X$ and any choice of basepoint, $\phi _{E}(%
\widetilde{X})$ is uniformly open in $X_{E}$.
\end{enumerate}
\end{theorem}

\begin{proof}
$1\Rightarrow 2$ follows from Theorem 45 in \cite{BPUU}. For $2\Rightarrow 3$%
, let $E$ be an entourage in $X$, $A:=\phi _{E}(\widetilde{X})$. Suppose
that $\left( [\alpha ]_{E},[\beta ]_{E}\right) \in E_{A}^{\ast }$ and $%
[\alpha ]_{E}\in A$. So there is some $([\gamma _{D}]_{D})\in \widetilde{X}$
such that $[\gamma _{E}]_{E}=[\alpha ]_{E}$ and we may write $\alpha
:=\{x_{0},...,x_{n},x\}$ and $\beta :=\{x_{0},...,x_{n},y\}$, with $(x,y)\in
E_{A}$. This in turn means that we have $([\alpha _{D}]_{D}),([\beta
_{D}]_{D})\in \widetilde{X}$ with endpoints $x$ and $y$ such that $([\alpha
_{E}]_{E}),([\beta _{E}]_{E})\in E^{\ast }\cap \left( A\times A\right) $. So
we may now write $\alpha _{E}:=\{y_{0},...,y_{m},x\}$ and $\beta
_{E}:=\{y_{0},...,y_{m},y\}$ and $(x,y)\in E$. Consider $([\gamma _{D}\ast
\alpha _{D}^{-1}\ast \beta _{D}]_{D})\in \widetilde{X}$. Using an $E$%
-homotopy like the one in the proof of Lemma \ref{form} we have%
\begin{equation*}
\lbrack \gamma _{E}\ast \alpha _{E}^{-1}\ast \beta
_{E}]_{E}=[x_{0},...,x_{n},x,y_{m},...,y_{0},y_{1},...,y_{m},y]_{E}
\end{equation*}%
\begin{equation*}
=[x_{0},...,x_{n},x,y]_{E}=[\beta ]_{E}
\end{equation*}%
This implies that $[\beta ]_{E}\in A$ and finishes the proof that $A$ is
uniformly open.

To prove $3\Rightarrow 1$, note that by Proposition \ref{up}, $\widetilde{X}$
has an invariant universal basis with respect to the isomorphic action of $%
\delta _{1}(X)$. Corollary \ref{greatcor} and Lemma \ref{greatlem} together
show that $\phi $ is a bi-uniformly continuous surjection, hence a quotient
with respect to this action (cf. Theorem 11, \cite{P}). Lemma \ref{61} now
finishes the proof.
\end{proof}

The next corollary is the statement Corollary 61 in \cite{BPUU}:

\begin{corollary}
\label{new}If $f:X\rightarrow Y$ is a bi-uniformly continuous surjection
where $X$ is universal and $Y$ is uniform then $Y$ is coverable.
\end{corollary}

\begin{proof}
According to Proposition 57 in \cite{BPUU} we have the lift $%
f_{L}:X\rightarrow \widetilde{Y}$ which satisfies $\phi \circ f_{L}=f$,
where $\phi :\widetilde{Y}\rightarrow Y$ is the projection. But then $\phi $
must be a uniformly continuous surjection. If $E$ is an entourage in $%
\widetilde{Y}$, then since $f$ is bi-uniformly continuous, $f(f_{L}^{-1}(E))$
is an entourage that is contained in $\phi (E)$. This proves that $\phi $ is
bi-uniformly continuous and hence $Y$ is coverable by Theorem \ref{main}.
\end{proof}

\begin{corollary}
If $X$ is coverable then $E$ is a covering entourage if and only if $X_{E}$
is chain connected.
\end{corollary}

\begin{proof}
If $E$ is a covering entourage then by definition $\phi _{E}:\widetilde{X}%
\rightarrow X_{E}$ is surjective. Since $\widetilde{X}$ is chain connected,
so is $X_{E}$. Conversely, if $X_{E}$ is chain connected then by Lemma \ref%
{uolem} and the third part of Theorem \ref{main}, $\phi _{E}$ must be
surjective.
\end{proof}

Note that the argument $3\Rightarrow 1$ in the proof of Theorem \ref{main}
is constructive; it actually provides a covering basis. We can now sort
through the steps to help identify this basis. The proof of Lemma \ref{61},
which is actually in \cite{BPUU}, shows that the covering entourages are of
the form $\phi (G)$, where $G$ is an invariant universal entourage in $%
\widetilde{X}$. The universal entourages in $\widetilde{X}$ come from
Proposition \ref{up}, and they are of the form $\phi _{E}^{-1}(E^{\ast })$
for any $E$. Letting $A:=\phi _{E}(\widetilde{X})$ we have 
\begin{equation*}
\phi (\phi _{E}^{-1}(E^{\ast }))=\phi (\phi _{E}^{-1}(E^{\ast }\cap (A\times
A))=\phi _{XE}\circ \phi _{E}\circ \phi _{E}^{-1}(E^{\ast }\cap \left(
A\times A\right) )=E_{A}
\end{equation*}%
Note that $E_{A}$ is chain connected since $\widetilde{X}$ is. Combining
this with Corollary \ref{cool} we obtain:

\begin{corollary}
\label{extract}Let $X$ be a coverable uniform space. For any entourage $E$, $%
X_{E}$ has a unique chain connected uniformly open set $A$ containing the
basepoint, and $E_{A}$ is a covering entourage.
\end{corollary}

\begin{example}
\label{exex}We will illustrate how Corollary \ref{extract} extracts a
covering entourage from a non-covering entourage in the topological group $%
\mathbb{R}$. In a topological group with left uniformity, entourages are
completely determined by symmetric open subsets of the identity (which
always serves as the basepoint). For example, in $\mathbb{R}$, if $U$ is any
such set, there corresponds an entourage $E(U):=\{(x,y):x-y\in U\}$. The
open set corresponding to $E(U)^{\ast }$ in $\mathbb{R}_{U}:=\mathbb{R}%
_{E(U)}$ is denoted by $U^{\ast }$. Using open sets rather than entourages
makes it easier to see what is going on. In Example 48 of \cite{BPTG} we
considered the set $U:=(-1,1)\cup (2,4)\cup (-4,-2)$. In this example the
components of $U$ are far enough apart that $\mathbb{R}_{U}$ consists of the
topological group $\mathbb{R\times Z}$ with the product uniform structure ($%
\mathbb{Z}$ is discrete). The idea here is that the two outer components of $%
U$ cannot be reached from $0$ by $(-1,1)$-chains, so for example the
equivalence class of the chain $\{0,3\}$ lies in a different component from
the identity component $A:=\mathbb{R}\times \{0\}$ in $\mathbb{R}_{U}$.
Along these lines, it is not hard to show that 
\begin{equation*}
U^{\ast }=(-1,1)\times \{0\}\cup (2,4)\times \{1\}\cup (-4,-2)\times \{-1\}
\end{equation*}%
That is, the two outer components of $U^{\ast }$ do not lie in $A$, which
clearly is the unique chain connected uniformly open set containing the
identity in $\mathbb{R\times Z}$. Now we have $\phi _{\mathbb{R}%
E(U)}(U^{\ast }\cap A)=(-1,1):=V$. Since $V$ is connected and $\mathbb{R}$
is simply connected, $E(V)$ is a covering entourage (cf. \cite{BPTG}).
\end{example}

In \cite{RC} the notion of Rips complex is extended from metric spaces to
uniform spaces: $R(X,E)$ is the subcomplex of the full complex over $X$
having as simplices all $\{x_{0},...,x_{n}\}$ such that $(x_{i},x_{j})\in E$
for all $i$ and $j$. According to (\cite{S}, Section 3.6), any path in $%
R(X,E)$ is, up to homotopy, uniquely identified with a simplicial path,
which in turn is uniquely determined by its vertices. These vertices,
obviously, form an $E$-chain, and the basic moves in a fixed-endpoint
simplicial homotopy of simplicial paths (adding or removing a pair of edges
that span $2$-simplex with one edge already in the path) correspond
precisely to the basic moves in an $E$-homotopy (adding or removing a point
so as to preserve that one has an $E$-chain). That is, the set of all
fixed-endpoint homotopy equivalence classes of paths in $R(X,E)$ starting at
a base point $\ast $ is naturally identified with $X_{E}$.

Using this natural identification of fixed-endpoing homotopies of paths in $%
R(X,E)$ and $E$-homotopies of $E$-chains, we will translate the basic
definitions of \cite{RC}. Two $E$-chains $\alpha :=\{x_{0},...,x_{n}\}$ and $%
\beta :=\{y_{0},...,y_{k}\}$ (maybe without the same endpoints) are said in 
\cite{RC} to be $E$-homotopic if (1) $(x_{0},y_{0}),(x_{n},y_{k})\in E$ and
(2) $\beta $ is (fixed-endpoint) $E$-homotopic to $%
\{y_{0},x_{0},...,x_{n},y_{k}\}$. If $\alpha $ and $\beta $ have the same
pair of endpoints then of course \textquotedblleft $E$-homotopic%
\textquotedblright\ has the same meaning as in \cite{BPUU}. If $\alpha $ and 
$\beta $ only have the same \textit{starting} point $x_{0}=y_{0}=\ast $,
then it is easy to check that $\alpha $ and $\beta $ are $E$-homotopic
precisely when $([\alpha ]_{E},[\beta ]_{E})\in E^{\ast }$. A \textit{%
generalized curve from }$x$ to $y$ is defined in \cite{RC} to be a
collection $\{[c_{E}]_{E}\}$ of $E$-homotopy classes of $E$-chains joining $%
x $ and $y$ such that if $F\subset E$ then $\left[ c_{F}\right]
_{E}=[c_{E}]_{E}$. The set of all generalized curves starting at $\ast $ is
called $GP(X,\ast )$ in \cite{RC}, but this set is obviously none other than 
$\widetilde{X}$ via the identification $\{[c_{E}]_{E}\}\leftrightarrow
\left( \lbrack c_{E}]_{E}\right) $. The authors define a \textquotedblleft
natural uniform structure\textquotedblright\ on $GP(X,\ast )$ (a.k.a. $%
\widetilde{X}$) by taking, for each entourage $F$ in $X$, the set of all
pairs $\left( \{[c_{E}]_{E}\},\{[d_{E}]_{E}\}\right) $ such that $c_{F}$ is $%
F$-homotopic to $d_{F}$. Since $c_{F}$ and $d_{F}$ both start at $\ast $, as
pointed out above this is equivalent to $([c_{F}]_{F},[d_{F}]_{F})\in
F^{\ast }$. That is, the basis that they define consists precisely of the
sets $\phi _{F}^{-1}(F^{\ast })$, which of course is a basis for the inverse
limit uniform structure on $\widetilde{X}$. In other words, $GP(X,\ast )$
and $\widetilde{X}$ are one and the same space. Moreover, the mapping $\pi
_{X}:GP(X,\ast )\rightarrow X$ of \cite{RC} is the endpoint mapping
(identical to $\phi :\widetilde{X}\rightarrow X$), and the uniform
fundamental group $\overset{\vee }{\pi }(X,\ast )$ of \cite{RC} is $\pi
_{X}^{-1}(\ast )=\phi ^{-1}(\ast )$ (\textquotedblleft generalized
loops\textquotedblright ) with operation induced by concatenation (identical
to $\delta _{1}(X)$).

According to \cite{RC}, a uniform space $X$ is called \textit{joinable}\ if
every pair of points in $X$ is joined by a generalized curve; clearly this
is equivalent to the surjectivity of $\phi :\widetilde{X}\rightarrow X$. In 
\cite{RC} $X$ is called \textit{uniformly joinable}\ if for every entourage $%
E$ there is an entourage $F$ such that whenever $(x,y)\in F$, $x$ and $y$
are joined by a generalized curve $\{[c_{D}]_{D}\}$ that is
\textquotedblleft $E$-short\textquotedblright\ in the sense that $%
[c_{E}]_{E}=[\{x,y\}]_{E}$.

\begin{corollary}
\label{eqrc}If $X$ is a chain connected uniform space then $X$ is coverable
if and only if $X$ is uniformly joinable.
\end{corollary}

\begin{proof}
Suppose that $\phi _{E}(\widetilde{X})$ is uniformly open for all $E$; so
there is some $F\subset E$ such that if $([\alpha ]_{E},[\beta ]_{E})\in
F^{\ast }$ and $[\alpha ]_{E}\in \phi _{E}(\widetilde{X})$ then $[\beta
]_{E}\in \phi _{E}(\widetilde{X})$. Let $(x,y)\in F$. Since $X$ is chain
connected there is some $F$-chain $\alpha =\{x_{0}=\ast ,...,x_{n-1},x\}$
and we may let $\beta :=\{x_{0},...,x_{n-1},x,y\}$. Note that since 
\begin{equation*}
([x_{0},...,x_{i}]_{E},[x_{0},...,x_{i},x_{i+1}]_{E})\in F^{\ast }
\end{equation*}%
for all $i$, it follows by induction on $i$ that $[\alpha ]_{E}\in \phi _{E}(%
\widetilde{X})$. Likewise $[\beta ]_{E}\in \phi _{E}(\widetilde{X})$ and we
have $[\alpha ]_{E}=[\alpha _{E}]_{E}$ for some $([\alpha _{D}]_{D})\in 
\widetilde{X}$ and $[\beta ]_{E}=[\beta _{E}]_{E}$ for some $([\beta
_{D}]_{D})\in \widetilde{X}$. But the concatenated generalized curve $%
\{[\alpha _{D}^{-1}\ast \beta _{D}]_{D}\}$ certainly satisfies the $E$-short
condition $[\alpha _{E}^{-1}\ast \beta _{E}]_{E}=[\{x,y\}]_{E}$; in fact one
may remove the points $x_{0},x_{1},x_{1},...x_{n-1},x_{n-1},x$ in succession
to create an $E$-homotopy between $\alpha _{E}^{-1}\ast \beta _{E}$ and $%
\{x,y\}$.

If $X$ is uniformly joinable and $E$ is an entourage, by definition there is
some entourage $F\subset E$ such that if $(x,y)\in F$, $x$ and $y$ are
joined by an $E$-short generalized curve. Let $\alpha =\{\ast
=x_{0},...,x_{n}\}$ be an $E$-chain with $[\alpha ]_{E}\in \phi _{E}(%
\widetilde{X})$. If $([\alpha ]_{E},[\beta ]_{E})\in F^{\ast }$ then by
definition of $F^{\ast }$ we may assume that $\beta $ is of the form $\{\ast
=x_{0},...,x_{n-1},x\}$ with $(x,x_{n})\in F$. That is, there is an $E$%
-short generalized curve $\{[c_{D}]_{D}\}$ joining $x_{n}$ and $x$ with $%
c_{E}=\{x_{n},x\}$. Now if $\phi _{E}([\alpha _{D}]_{D})=[\alpha ]_{E}$ then 
$g:=([\alpha _{D}\ast c_{D}]_{D})\in \widetilde{X}$ satisfies $\phi
_{E}(g)=[\beta ]_{E}$.
\end{proof}

\begin{remark}
In light of Corollaries \ref{greatcor} and \ref{eqrc} we have a very nice
way to distinguish between \textit{joinable and uniformly joinable for a
chain connected uniform space }$X$, namely that for a joinable space, $\phi
_{E}$ has uniformly open image for some $E$, while for a uniformly joinable
space, $\phi _{E}$ has uniformly open image for all $E$.
\end{remark}

\begin{remark}
Note that the equivalence of Theorem \ref{main}.2 and uniform joinability
for a chain connected space was proved in \cite{RC} using completely
different arguments.
\end{remark}


\begin{thebibliography}{9}
\bibitem{BPTG} V. Berestovskii and C. Plaut, Covering group theory for
topological groups, Top. Appl. 114 (2001) 141-186.

\bibitem{BPLC} V. Berestovskii and C. Plaut, Covering group theory for
locally compact groups, Top. Appl. 114 (2001) 187-199.

\bibitem{BPUU} V. Berestovskii and C. Plaut, Uniform Universal Covers of
Uniform Spaces, Top. Appl. 154 (2007) 1748--1777.

\bibitem{BPCR} V. Berestovskii and C. Plaut, Covering $\mathbb{R}$-trees,
preprint arXiv:0707.3609.

\bibitem{RC} N. Brodskiy, J. Dydak, B. Labuz, and A. Mitra, Rips Complexes
and universal covers in the uniform category, preprint arXiv:0706.3937.

\bibitem{P} C. Plaut, Quotients of uniform spaces, Top. Appl. 153 (2006)
2430-2444.

\bibitem{PL} C. Plaut, Lectures on Uniform Universal Covers of Uniform
Spaces (joint work with V. Berestovskii), University of Tennessee, Feb.,
2007.

\bibitem{S} E. Spanier, \textit{Algebraic Topology}, McGraw-Hill, New York,
1966.
\end{thebibliography}
\end{document}